\documentclass[11pt]{article}

\title{Random graph's Hamiltonicity is strongly tied to its minimum degree }

\author{
Yahav Alon
\thanks{‡School of Mathematical Sciences, Raymond and Beverly Sackler Faculty of Exact Sciences, Tel Aviv University,
Tel Aviv, 6997801, Israel. Email: yahavalo@mail.tau.ac.il.}
\and Michael Krivelevich
\thanks{School of Mathematical Sciences, Raymond and Beverly
Sackler Faculty of Exact Sciences, Tel Aviv University, Tel Aviv,
6997801, Israel. Email: krivelev@post.tau.ac.il. Partially supported by USA-Israel BSF grant 2014361, and by ISF grant 1261/17.}
}

\usepackage{amsmath,amsthm, amssymb,latexsym, amsfonts}

\begin{document}
\maketitle
\newtheorem{thm}{Theorem}
\newtheorem{propos}{Proposition}
\newtheorem{defin}{Definition}
\newtheorem{lemma}{Lemma}[section]
\newtheorem{corol}{Corollary}
\newtheorem{thmtool}{Theorem}[section]
\newtheorem{corollary}[thmtool]{Corollary}
\newtheorem{lem}[thmtool]{Lemma}
\newtheorem{defi}[thmtool]{Definition}
\newtheorem{prop}[thmtool]{Proposition}
\newtheorem{clm}[thmtool]{Claim}
\newtheorem{conjecture}{Conjecture}
\newtheorem{problem}{Problem}
\newcommand{\Proof}{\noindent{\bf Proof.}\ \ }
\newcommand{\Remarks}{\noindent{\bf Remarks:}\ \ }
\newcommand{\Remark}{\noindent{\bf Remark:}\ \ }

\begin{abstract}
We show that the probability that a random graph $G\sim G(n,p)$ contains no Hamilton cycle is $(1+o(1))Pr(\delta (G) < 2)$ for all values of $p = p(n)$. We also prove an analogous result for perfect matchings.
\end{abstract}

\section{Introduction and main results}\label{sec-intro}

Hamilton cycles are a central topic in modern graph theory, a fact that extends to the field of random graphs as well, with numerous and diverse results regarding the appearance of Hamilton cycles in random graphs obtained over recent years.\\
A classical result by Koml{\'o}s and Szemer{\'e}di \cite{KS83}, and independently by Bollob{\'a}s \cite{B84}, states that a random graph $G\sim G(n,p)$, with $np-\ln n - \ln \ln n \rightarrow \infty$, is asymptotically almost surely  Hamiltonian. It should also be noted that if $np -\ln n - \ln \ln n \rightarrow -\infty$ then asymptotically almost surely $\delta (G) \leq 1$, and thus $G$ is not Hamiltonian.\\
The same exact statement is true if one replaces the graph property \textit{$G$ is Hamiltonian} with the property $\mathit{\delta (G) \geq 2}$. This indicates a possible connection between the two properties, an indication made explicit when considering a stronger result proved by Bollob{\'a}s in \cite{B84} and by Ajtai, Koml{\'o}s and Szemer{\'e}di in \cite{AKS85}, regarding the \textit{hitting time} of Hamiltonicity.\\
Consider a \textit{random graph process}, defined as a sequence of (random) graphs on $n$ vertices $\tilde{G}(\sigma ) = \lbrace G_i(\sigma ) \rbrace _{i=0} ^{\binom{n}{2}}$, where $\sigma$ is an ordering of the edges of $K_n$ chosen randomly and uniformly from among all $\binom{n}{2}!$ such orderings. Set $G_0 (\sigma )$ to be a graph with no edges, and for all $1\leq i \leq \binom{n}{2}$, $G_i (\sigma )$ is obtained by adding the $i$--th edge according to the order $\sigma$ to  $G_{i-1} (\sigma )$. The \textit{hitting time} of a monotone and non--empty graph property $\mathcal{P}$, which we will denote as $\tau _{\mathcal{P}}( \tilde{G} (\sigma ))$, is a random variable equal to the index $i$ for which $G_i(\sigma ) \in \mathcal{P}$ but $G_{i-1}(\sigma ) \notin \mathcal{P}$.\\
Denote by $\mathcal{H}$ the property of Hamiltonicity, and by $\mathcal{D}_2$ the property of having minimum degree at least two. The result states that asymptotically almost surely $\tau _{\mathcal{D}_2}(\tilde{G}(\sigma )) = \tau _{\mathcal{H}}(\tilde{G}(\sigma ))$.\\
Indeed, the result about Hamiltonicity in $G(n,p)$ can be derived directly from the hitting time result, thus making the hitting time result stronger. In addition, the hitting time result indicates an explicit connection between the minimum degree of a random graph and the existence of a Hamilton cycle. The random graph process asymptotically almost surely becomes Hamiltonian at the exact same moment it has minimum degree at least two. This shows that minimum degree less than two is typically the chief obstacle for Hamiltonicity.\\
In light of this, it seems natural to ask whether the connection between the two properties can be expressed more explicitly. Furthermore, since the latter result revolves around the hitting time, one can ask whether this connection can also be observed in random graphs that are much more dense than the threshold density of $\frac{1}{2}(\ln n + \ln \ln n + \omega (n))$.\\
A partial answer to these questions has been given by McDiarmid and Yolov in 2016 \cite{MY16}. Based on a result by Hefetz, Krivelevich and Szab{\'o} \cite{HKS09}, the authors proved the following:\\
If $p\leq \frac{1}{2}$ is such that $\frac{p\cdot n \ln \ln \ln \ln n}{\ln n \cdot \ln \ln \ln n} \rightarrow \infty$, then the probability of $G\sim G(n,p)$ failing to contain a Hamilton cycle is at most $(1-p)^n\cdot \exp \left( O \left( \frac{\ln n \ln \ln \ln n}{\ln \ln \ln \ln n} \right) \right) $.\\
One can observe that the probability of $G \sim G(n,p)$ having $\delta (G) <2$ is of order $\Theta \left( np(1-p)^n \right) $, and so this result gives an explicit bound on the ratio between the probabilities of the negations of these two properties:\\
$$
1 \leq \frac{Pr(G\ \mbox{is\ not\ Hamiltonian})}{Pr(\delta (G) <2)} \leq \exp \left( O \left( \frac{\ln n \ln \ln \ln n}{\ln \ln \ln \ln n} \right) \right)
$$
This however only partially answers our question due to two gaps: first, it does not cover all of our target range of $p(n)$; second, it is far from being tight.\\
Here we close both these gaps, by proving the following main result:

\begin{thm}\label{main-thm}
Let $0 \leq p=p(n) \leq 1$, and let $G \sim G(n,p)$. Then $Pr(G\ \emph{is\ not\ Hamiltonian})=(1+o(1))Pr(\delta (G) < 2)$.
\end{thm}
\noindent By proving Theorem \ref{main-thm} we cover all possible values of $p$, as well as achieve the asymptotically tight ratio of $1+o(1)$.\\
Going forward, we utilize some of the approaches we use to prove Theorem \ref{main-thm} in order to sketch a proof of an analogous result regarding perfect matchings in a random graph.\\
Assume for simplicity that $n$ is even. Similar to the result by Koml{\'o}s and Szemer{\'e}di and of Bollob{\'a}s about the threshold probability of Hamiltonicity, a very early work by Erd\H{o}s and R{\'e}nyi \cite{ER} shows that whenever $np- \ln n \rightarrow \infty$, the random graph $G \sim G(n,p)$ asymptotically almost surely contains a perfect matching. Similar to the connection between Hamiltonicity and minimum degree 2, this statement is also true when replacing the property of containing a perfect matching with that of containing no isolated vertices.\\
In fact, it was later proved by Bollob{\'a}s and Thomason \cite{BT} that the hitting times of the two properties are asymptotically almost surely exactly equal to each other.\\
Utilizing some similar approaches to those used in the proof of Theorem \ref{main-thm} we sketch a proof the following:

\begin{thm}\label{matching-thm}
Let $0 \leq p=p(n) \leq 1$, and let $G\sim G(n,p)$. Then $Pr(G\ \emph{contains\ no\ perfect\ matching})=(1+o(1))Pr(\delta (G) = 0)$.
\end{thm}
\noindent Once more this theorem covers the full range and provides a very tight ratio.\\

This paper is structured as follows: in Section \ref{sec-defs} we list the notations and definitions to be used throughout the paper, as well as some auxiliary results needed for our proofs. In Section \ref{sec-proof} we provide a proof of Theorem \ref{main-thm}. In Section \ref{sec-PM} we sketch a proof of Theorem \ref{matching-thm}.\\
Portions of the proofs in Section \ref{sec-proof} and Section \ref{sec-PM} and several of the techniques we employ on our proofs are inspired by \cite{K}.

\section{Notation, definitions and auxiliary results}\label{sec-defs}

In this section we provide several definitions and results to be used in the following sections.\\
Throughout the paper, it is assumed that all logarithmic functions are in the natural base, unless explicitly stated otherwise.\\
We suppress the rounding notation occasionally to simplify the presentation.\\
The following standard graph theoretic notations will be used:
\begin{itemize}
\item $N_G(U)$ : the external neighbourhood of a vertex subset $U$ in the graph $G$, i.e.
$$
N_G(U) = \lbrace v \in V(G)\setminus U:\ v\ \mbox{has\ a\ neighbour\ in}\ U \rbrace.
$$
\item $e_G(U)$: the number of edges spanned in a graph $G$ by a vertex subset $U$. This will sometimes be abbreviated as $e(U)$, when the identity of $G$ is clear from the context.
\item $e_G(U,W)$: the number of edges of $G$ between the two disjoint vertex sets $U,W$. This will sometimes be abbreviated as $e(U,W)$ when $G$ is clear from the context.
\item $E_G(v)$: the set of edges in a graph $G$ incident to the vertex $v$.
\end{itemize}

\subsection{Graph theory}\label{subsec-graphs}

\begin{defi}\label{expander}
Let $G=(V,E)$ be a graph, and let $\alpha>0$ and $k$ a positive integer. The graph $G$ is a {\em $(k,\alpha )$-expander} if $|N_G(W)|\geq \alpha |W|$ for every vertex subset $W\subset V$, $|W|\leq k$.
\end{defi}

\begin{defi}\label{def-booster}
Let $G=(V,E)$ be a graph. A non-edge $(u,v)\in E(G)$ is called a {\em booster} if the graph $G^{\prime}$ with edge set $E(G^{\prime})=E(G)\cup \lbrace (u,v) \rbrace$ is either Hamiltonian or has a path longer than a longest path of $G$.
\end{defi}

\begin{lem}\label{lemma-booster}
\emph{(P{\'o}sa 1976 \cite{POS})} Let $G$ be a connected non--Hamiltonian graph, and assume that $G$ is a $(k,2)$--expander. Then $G$ has at least $\frac{(k+1)^2}{2}$ boosters.
\end{lem}

\begin{defi}\label{def-hamconn}
A graph $G$ is \emph{Hamilton--connected} if for every two vertices $u,v\in V(G)$, $G$ contains a Hamilton path with $u,v$ as its two endpoints.
\end{defi}

\begin{thmtool}\label{chvatal-erdos}
{\em (Chv{\'a}tal--Erd\H{o}s Theorem \cite{CE})} Let $G=(V,E)$ be a graph such that $\alpha (G) < \kappa (G)$. Then $G$ is Hamilton--connected.
\end{thmtool}

\subsection{Binomial coefficients and binomial distribution}\label{subsec-binom}

\begin{lem}\label{coefficients}
Let $1\leq l \leq k \leq n$ be integers. Then the following inequalities hold:
\begin{enumerate}
\item $\binom{n}{k} \leq \left( \frac{en}{k} \right) ^k$\,;
\item $\frac{\binom{n-l}{k}}{\binom{n}{k}} \leq e^{-\frac{l\cdot k}{n}}$\,.
\end{enumerate}
\end{lem}

\begin{lem}\label{binom-rv}
Let $1\leq k \leq n$ be integers, $0 < p < 1$, and let $X\sim Bin(n,p)$. Then the following inequalities hold:
\begin{enumerate}
\item $Pr(X \geq k) \leq \left( \frac{enp}{k} \right) ^k$\,;
\item $Pr(X = k) \leq \left( \frac{enp}{k(1-p)} \right) ^k\cdot e^{-np}$\,.
\end{enumerate}
\end{lem}

\begin{lem}\label{chernoff}
{\em (Chernoff bound for binomial tails, see e.g. \cite{CHER})} Let $X\sim Bin(n,p)$. Then the following inequalities hold:
\begin{enumerate}
\item $Pr(X <(1-\delta )np) \leq \exp \left( -\frac{\delta ^2 np}{2} \right)$ for every $\delta > 0$\,;
\item $Pr(X >(1+\delta )np) \leq \exp \left( -\frac{\delta ^2 np}{3} \right)$ for every $0 < \delta < 1$\,;
\item $Pr(X >(1+\delta )np) \leq \exp \left( -\frac{\delta np}{3} \right)$ for every $\delta > 0$\,.
\end{enumerate}
\end{lem}

\subsection{General auxiliary results}\label{subsec-general}

\begin{lem}\label{bonferroni}
{\em (Particular case of Bonferroni inequality)} Let ${\lbrace A_k \rbrace}_{k=1} ^n$ be a family of events in a probability space. Then:
$$
Pr\left( \bigcup_{k=1}^n A_k \right) \geq \sum_{k=1}^n Pr(A_k)\ - \sum_{1\leq k<l \leq n} Pr(A_k \cap A_l).
$$
\end{lem}

\section{Proof of main theorem}\label{sec-proof}

In this section we provide a proof for Theorem \ref{main-thm}. The theorem's statement covers all possible values of the edge probability $p(n)$, with the corresponding random graphs $G(n,p)$ having somewhat different characteristics in different parts of this range.\\
Since this is the case, we divide our proof into four parts, each corresponding to a different range of $p$. The general approach in each of the proofs, bar the first part, will be as thus: we will present a finite set of graph properties, say, $\lbrace \mathcal{P}_i \rbrace _{i\geq 0}$, with $\mathcal{P}_0$ being the property of having minimum degree at least 2, and show that:
\begin{itemize}
\item[(i)] For a random graph $G\sim G(n,p)$ and for all $i>0$, the probability of $G\notin \mathcal{P}_i$ is $o\left( Pr(G \notin \mathcal{P}_0) \right)$;
\item[(ii)] Any graph $G$ such that $G\in \bigcap_{i\geq 0} \mathcal{P}_i$ is Hamiltonian.
\end{itemize}
Since a Hamiltonian graph $G$ must have a minimum degree at least 2, meaning $Pr(\delta (G) < 2) \leq Pr(G\ \mbox{is\ not\ Hamiltonian})$, combining the two claims yields the theorem.\\
The four parts of the proof will be the following:

\begin{itemize}
\item {\bf The very sparse case}: let $p = p(n)$ be such that $np - \log n - \log \log n$ does not tend to infinity. In this case there is nothing to prove, since the theorem's statement is already known to be true. This is due to the result by Koml{\'o}s and Szemer{\'e}di, stating that for $G \sim G(n,p)$:
\[
\lim _{n \rightarrow \infty }\Pr(G\ \mbox{is Hamiltonian})=
\begin{cases}
    0 & np - \log n - \log \log n \rightarrow -\infty,\\
    e^{-e^{-c}} & np - \log n - \log \log n \rightarrow c,\\
    1 & np - \log n - \log \log n \rightarrow \infty.
\end{cases}
\]
It follows that $\lim _{n \rightarrow \infty }\Pr(G\ \mbox{is not Hamiltonian}) = \lim _{n \rightarrow \infty }\Pr(\delta (G) < 2) = C > 0$, and the statement therefore holds.
\item {\bf The sparse case}: in Section \ref{subsec-sparse} we provide the proof of the theorem, assuming $np- \log n - \log \log n \rightarrow \infty$ and $p \leq \frac{100\log n}{n}$. In the proof of this case, we aim to show that with the appropriate probability for $G\sim G(n,p)$, $G$ contains an $\left( \frac{n}{4},2 \right)$--expander as a subgraph that is sparse, relative to the average degree of $G$. Then invoking Lemma \ref{lemma-booster} we conclude that this subgraph can be made Hamiltonian with high probability by adding boosters that are also edges of $G$.\\
In this range of $p$, our main concern is the existence of many vertices of degree sublinear in $np$. This concern is addressed by showing that with high probability these vertices are relatively few and far apart.\\
This proof, up to the use of the $G(n,p)$ model and the tightness of the bounds sought, is essentially due to \cite{K}.
\item {\bf The dense case}: in Section \ref{subsec-dense} we cover the case $\frac{100\log n}{n} \leq p \leq 0.01$. The proof of this case is fairly similar to the proof of the sparse case, and in fact somewhat simpler due to the fact that it is now highly unlikely for $G$ to contain more than one vertex of degree sublinear in $np$.
\item {\bf The very dense case}: in Section \ref{subsec-vdense} we provide the proof of the theorem, assuming $p \geq 0.01$. In this range the graph is indeed typically very dense, and so removing a minimum degree vertex will typically result in a graph whose independence number is smaller than its vertex connectivity. This fact allows for the use of the Chv{\'a}tal-Erd\H{o}s Theorem to prove Hamiltonicity of the original graph.
\end{itemize}

\subsection{The sparse case}\label{subsec-sparse}

Recall that we assume here that $np - \log n - \log \log n \to \infty$ and $p \leq \frac{100\log n}{n}$. Define
$$
d_0 = 0.001\log n,
$$
and for a graph $G=(V,E)$, denote
$$
SMALL(G) = \lbrace v\in V(G):\, d(v) \leq d_0 \rbrace.
$$

Now, define the following graph properties:

\begin{itemize}
\item[{\bf (P0)}] $\delta (G) \geq 2$;
\item[{\bf (P1)}] $\Delta (G)\leq 800\log n$;
\item[{\bf (P2)}] $|SMALL(G)|\leq n^{0.3}$;
\item[{\bf (P3)}] $SMALL(G)$ does not contain two vertices $u,v$ such that $dist(u,v)\leq 4$;
\item[{\bf (P4)}] $\forall U \subset V(G)\ s.t.\ |U|\leq \frac{n}{\sqrt{\log n}}:\ e(U) < |U|{\log }^{\frac{3}{4}} n$;
\item[{\bf (P5)}] $\forall U,W\subset V(G)\ disjoint\ s.t.\ |U|\leq \frac{n}{\sqrt{\log n}},\ |W|\leq |U|{\log }^{\frac{1}{4}} n:\ e(U,W)<\frac{d_0|U|}{2}$;
\item[{\bf (P6)}] $\forall U,W\subset V(G)\ disjoint\ s.t.\ |U|=|W|= \frac{n}{\sqrt{\log n}} : e(U,W)\geq \frac{1}{2}n$;
\item[{\bf (P7)}] Every $\left( \frac{n}{4},2 \right) $-expanding subgraph $\Gamma$ of $G$ with at most $(d_0 + 1)n$ edges is either Hamiltonian, or has a booster in $E(G) \setminus E(\Gamma )$.
\end{itemize}

\begin{lemma}\label{sparse-properties}
Let $p=p(n)=\frac{\log n + \log \log n + \omega (n)}{n}$, where $\omega (n)$ is such that $\lim _{n \rightarrow \infty}\omega (n) = \infty$ and $p\leq \frac{100\log n}{n}$, and let $G \sim G(n,p)$ be a random graph. Then the probability that all properties {\bf (P0)}--{\bf (P7)} hold is $1-(1+o(1))Pr(G \notin {\bf{(P0)}})$.
\end{lemma}

\noindent {\bf Proof.} Clearly, the probability that exists a property among {\bf (P0)}--{\bf (P7)} that does not hold is at least $Pr(G\notin {\bf (P0)})$, and at most $\sum_{i=0}^7Pr(G\notin {\bf (Pi)})$. We will prove that for all $1\leq i \leq 7$ one has $Pr(G\notin {\bf (Pi)}) = o(Pr(G\notin {\bf (P0)}))$, and thus will establish the lemma.

First, we bound $Pr(G\notin {\bf (P0)})$ from below:

\noindent For a vertex $v\in V(G)$, denote by $A_v$ the event ``$d(v)<2$". So $Pr(G\notin {\bf (P0)}) = Pr\left( \bigcup_{v\in V(G)}A_v \right)$. By Bonferroni's inequality (Lemma \ref{bonferroni}):
$$
Pr(G\notin {\bf (P0)}) \geq \sum_{v\in V(G)}Pr(A_v)\ -\sum_{u,v\in V(G)} Pr(A_u \cap A_v).
$$
Let $u,v\in V(G)$. We first observe that
\[
\begin{array}{rcl}
Pr(A_v) & \geq & Pr(d(v)=1) \\
& = & (n-1)p(1-p)^{n-2} \\
& \geq & e^{-\frac{p}{1-p}(n-2)}\cdot (n-1)p \\
& \geq & e^{-np+O(np^2)}np.\\
\end{array}
\]

Since $\log n < np \leq 100\log n$, we get:
$$
Pr(A_v)  =  e^{-np+O(np^2)}np \geq \frac{1}{n}\cdot e^{-\omega (n) + o(1)}.
$$

Now, the probability of $A_u\cap A_v$ is at most the probability that $e_G(\lbrace u,v \rbrace ,V\setminus \lbrace u,v \rbrace ) \leq 2$:

\[
\begin{array}{rcl}
Pr(A_u\cap A_v) & \leq & Pr(e_G(\lbrace u,v \rbrace ,V\setminus \lbrace u,v \rbrace ) \leq 2) \\
& = & (1-p)^{2n-6}\left( (1-p)^2 + (2n-4)p(1-p) + \binom{2n-4}{2}p^2 \right) \\
& \leq & e^{-2np + o(1)}\cdot O(n^2p^2)\\
& = & O\left( \frac{1}{n^2}\cdot e^{-2\omega (n)} \right).
\end{array}
\]

And so overall we get:
$$
Pr(G\notin {\bf (P0)}) \geq e^{-\omega (n) + o(1)}-\binom{n}{2}\cdot O\left( \frac{1}{n^2}\cdot e^{-2\omega (n)} \right) = (1-o(1))e^{-\omega (n)}.
$$

For the rest of the properties we bound their probabilities from above:

\begin{itemize}
\item[{\bf (P1).}] Let $v\in V(G)$. By the union bound:
$$
Pr(G\notin {\bf (P1)}) \leq n\cdot Pr(d(v)\geq 800\log n),
$$
and by Lemma \ref{binom-rv} this gives:
$$
Pr(G\notin {\bf (P1)}) \leq n\left( \frac{enp}{800\log n} \right)^{800\log n} \leq n^{-799} = o(Pr(G\notin {\bf (P0)})).
$$

\item[{\bf (P2).}] The probability of $|SMALL(G)| \geq n^{0.3}$ is less than the probability that exists a set $S$ of size $n^{0.3}$ with $e(S,V \setminus S) \leq d_0\cdot n^{0.3}$, so by the union bound:
\begin{eqnarray*}
Pr(G\notin {\bf (P2)}) & \leq & \binom{n}{n^{0.3}}\sum_{k=0}^{d_0n^{0.3}}Pr\left( e(S,V \setminus S)=k \right) \\
& \leq & \binom{n}{n^{0.3}}\left( d_0n^{0.3}+1\right) Pr\left( e(S,V \setminus S) =d_0n^{0.3}\right).
\end{eqnarray*}

By Lemma \ref{binom-rv} we can estimate the latter expression from above by:
\begin{eqnarray*}
n^{(0.7+o(1))n^{0.3}}{\left( \frac{en^{1.3}p}{d_0n^{0.3}(1-p)} \right) }^{d_0n^{0.3}} e^{-(1-o(1))n^{1.3}p}
\leq {\left( n^{-0.2} \left( \frac{100e}{0.001\cdot 0.9} \right) ^{0.001\log n} \right) }^{n^{0.3}}.\\
\end{eqnarray*}

Finally, since $\left( \frac{100e}{0.001\cdot 0.9} \right) ^{0.001} \leq e^{0.015}$, and since $\omega (n) = o(n^{0.3}\log n)$, we obtain:

$$
Pr(G\notin {\bf (P2)}) \leq n^{-0.1n^{0.3}} = o(Pr(G\notin {\bf (P0)})).
$$

\item[{\bf (P3).}] The probability that exist $u,v \in V(G)$ such that $u,v \in SMALL(G),\ dist(u,v) = k$ is at most the probability that there is a path $\mathcal{P}$ of length $k$ between them. By the union bound this is at most
$$
\binom{n}{k+1}p^k\cdot Pr(\mbox{Both endpoints of } \mathcal{P} \mbox{ have degree at most } d_0 \mid \mathcal{P} \in G).
$$
The probability that two endpoints $u,v$ of a path both have at most $d_0$ neighbours is at most the probability that they have at most $2(d_0-2)$ edges that are not $(u,v)$ or a part if the path. By Lemma \ref{binom-rv} we have
$$ Pr(u,v\in SMALL(G)) \leq (2d_0-3)\left( \frac{enp}{(d_0-2)(1-p)} \right) ^{2d_0-4}e^{-2p(n-6)} \leq n^{0.03}e^{-2np}.
$$

Apply the union bound once more, and sum over $1\leq k\leq 4$:
\begin{eqnarray*}
Pr(G\notin {\bf (P3)}) & \leq & \sum_{k=1}^4 \binom{n}{k+1}p^kn^{0.03}e^{-2np}\\
& \leq & \sum_{k=1}^4 n^{0.03}e^{-2np}\cdot n \cdot (100\log n)^k \\
& \leq & \frac{4\cdot 10^8\log ^4 n}{n^{0.9}}e^{-2\omega (n)} = o(Pr(G\notin {\bf (P0)})).
\end{eqnarray*}

\item[{\bf (P4).}] By the union bound and Lemma \ref{binom-rv}, the probability that exists a vertex set $U \subseteq V(G)$ of size $|U|=k\leq \frac{n}{\sqrt{\log n}}$ that contradicts {\bf (P4)}, namely $e(U) \geq k\log ^{\frac{3}{4}} n$, is at most
\begin{eqnarray*}
&& \binom{n}{k} \left( \frac{e\binom{k}{2}p}{k\log ^{\frac{3}{4}}n} \right)^{k\log ^{\frac{3}{4}}n}
\leq \left( \frac{en}{k}\cdot \left( \frac{50ek \cdot \log ^{\frac{1}{4}}n}{n} \right)^{\log ^{\frac{3}{4}}n} \right) ^k
\leq \left( \frac{50e}{\log ^{0.24}n} \right)^{k\log ^{\frac{3}{4}}n}\\
& \leq & (\log n)^{-0.2k\log ^\frac{3}{4} n}.
\end{eqnarray*}

Clearly for $U$ of size $k$ to have $e(U) \geq |U|\log ^{\frac{3}{4}} n$, it has to be that $k \geq \log ^{\frac{3}{4}}n$. So by the union bound:

$$
Pr(G\notin {\bf (P4)}) \leq \sum_{k=\log ^{\frac{3}{4}}n}^{\frac{n}{\sqrt{\log n}}} (\log n)^{-0.2k\log ^\frac{3}{4} n}
= (1+o(1))(\log n)^{-0.2\log ^{\frac{3}{2}}n},
$$

which is sufficiently small.

\item[{\bf (P5).}] Let $U,W$ be two distinct vertex subsets, with $|U| = k \leq \frac{n}{\sqrt{\log n}}$. Clearly, if $U,W$ contradict {\bf (P5)} and $|W| < k\log ^{\frac{1}{4}} n$ then there also exists $W^{\prime }$ such that $U,W^{\prime }$ contradict {\bf (P5)}, and $|W^{\prime }| = k\log ^{\frac{1}{4}} n$; so it can be assumed that $|W| = k\log ^{\frac{1}{4}} n$. By Lemma \ref{binom-rv} the probability that exist such $U,W$ with $e(U,W)\geq \frac{d_0k}{2}$ is at most
\begin{eqnarray*}
&& \binom{n}{k}\binom{n}{k\log ^{\frac{1}{4}} n}\left( \frac{2ek\log ^{\frac{1}{4}}n\cdot p}{d_0} \right)^{\frac{d_0k}{2}}
\leq \left( \frac{en}{k\log ^{\frac{1}{8}}n} \right)^{2k\log ^{\frac{1}{4}}n}\cdot \left( \frac{2000ekp}{\log ^{\frac{3}{4}}n} \right)^{0.0005k\log n} \\
& \leq & \left( \frac{k}{n} \right)^{0.0004 k \log n}\cdot \log ^{0.00015 k \log n} n = (\log n)^{-\Omega (k\log n)},
\end{eqnarray*}

and by the union bound we obtain
$$
Pr(G\notin {\bf (P5)}) \leq \sum_{k=1}^{\frac{n}{\sqrt{\log n}}} (\log n)^{-\Omega (k\log n)} = (1+o(1))(\log n)^{-\Omega (\log n)} = o(Pr(G\notin {\bf (P0)})).
$$

\item[{\bf (P6).}] By the union bound:
\begin{eqnarray*}
Pr(G \notin {\bf (P6)}) & = & \binom{n}{\frac{n}{\sqrt{\log m}}}^2 Pr\left( Bin \left( \frac{n^2}{\log n},p \right) < \frac{1}{2}n \right) \\
& \leq & e^{o(n)} Pr\left( Bin \left( \frac{n^2}{\log n},\frac{\log n}{n} \right) < \frac{1}{2}n \right). \\
\end{eqnarray*}

By applying Lemma \ref{chernoff} we bound  the expression above by:
$$
\exp \left( o(n) - \frac{n}{8} \right) = o(Pr(G\notin {\bf (P0)})).
$$

\item[{\bf (P7).}] By Lemma \ref{lemma-booster}, a non-Hamiltonian $\left( \frac{n}{4},2 \right)$-expander has at least $\frac{n^2}{32}$ boosters. By the union bound, the probability of the existence of such a subgraph $\Gamma\subseteq G$, with none of its boosters being an edge of $G$, is at most
\begin{eqnarray*}
Pr(G \notin {\bf (P7)}) & \leq & \sum_{k=1}^{(d_0+1)n} \binom{\binom{n}{2}}{k}p^k (1-p)^{\frac{n^2}{32}}\\
& \leq & O(n\log n)\cdot \left( \frac{enp}{2d_0} \right)^{(d_0+1)n} e^{-\frac{n^2p}{32}} \\
& \leq & O(n\log n)\cdot \left( \frac{100e}{0.002} \right)^{0.002n\log n} e^{-\frac{n\log n}{32}}\\
& \leq & \exp ((0.025 - 0.03)n\log n) = o(Pr(G\notin {\bf (P0)})).
\end{eqnarray*}
\end{itemize}
\hfill$\Box$

\begin{lemma}\label{sparse-P6prime}
Let $G$ be a graph such that {\bf(P0)}, {\bf(P1)} and {\bf(P6)} hold.  Then $G$ contains a subgraph $\Gamma_0$ such that:
\begin{enumerate}
\item $\delta (\Gamma_0) \geq 2$,
\item $\Gamma_0$ has at most $d_0n$ edges,
\item $\forall v\notin SMALL(G)\,:\, d_{\Gamma _0}(v) \geq d_0$,
\item $\forall U,W\subset V(G)\ disjoint\ vertex\ sets\ s.t.\ |U|=|W|= \frac{n}{\sqrt{\log n}} : e_{\Gamma _0}(U,W)\geq 1$.
\end{enumerate}
\end{lemma}

\noindent {\bf Proof.} Consider the following construction of a random subgraph $\Gamma _0$ of $G$ with at most $d_0n$ edges and minimum degree at least $\min \lbrace \delta (G), d_0 \rbrace $:\\
For each $v\in V(G)$ define $E_v = E_G(v)$ if $v \in SMALL(G)$, and let $E_v$ be a random subset of $E_G(v)$ of size exactly $d_0$ chosen uniformly if $v \notin SMALL(G)$. Finally, set $E(\Gamma _0) = \bigcup_{v\in V(G)}E_v$. Clearly, $\Gamma _0$ has at most $d_0n$ edges and minimum degree at least $\min \lbrace \delta (G), d_0 \rbrace $, and every vertex not in $SMALL(G)$ has degree at least $d_0$.\\
We will show that with positive probability (and in fact, with high probability) a random subgraph $\Gamma _0$ constructed in this manner is such that for every pair of disjoint sets $U,W\subseteq V(G)$, with $|U|=|W|= \frac{n}{\sqrt{\log n}}$, the graph $\Gamma _0$ has an edge between $U$ and $W$, thus ensuring the existence of the requested subgraph. By the union bound, the probability that exist such $U,W$ with no edge between them is at most
\begin{eqnarray*}
\binom{n}{\frac{n}{\sqrt{\log n}}}^2\prod_{u\in U}Pr(e_{\Gamma _0}(u,W)=0) & \leq & e^{o(n)}\cdot \prod _{u \in U} \frac{\binom{d_G(u) - e_G(u,W)}{d_0}}{\binom{d_G(u)}{d_0}} \\
& \leq &  e^{o(n)} \cdot \prod _{u \in U} e^{-\frac{d_0\cdot e_G(u,W)}{d_G(u)}}\\
& \leq & \exp \left( {-\frac{d_0}{\max _{u\in U}d_G(u)}\cdot e_G(U,W)+ o(n)} \right),
\end{eqnarray*}
the second inequality holding by Lemma \ref{coefficients}.\\
By {\bf (P1)} we have that $\max_{u\in U}d_G(u) \leq 800\log n$, and by {\bf (P6)} we have $e_G(U,W) \geq \frac{1}{2}n$, so the above expression is at most
$$
\exp \left( -\frac{0.001}{1600}n + o(n) \right) = o(1).
$$
\hfill$\Box$

\begin{lemma}\label{sparse-expander}
Let $G$ be a graph such that {\bf(P0)}--{\bf(P6)} hold.  Then $G$ contains a subgraph $\Gamma_0$ such that:
\begin{enumerate}
\item $\Gamma_0$ has at most $d_0n$ edges,
\item $\Gamma_0$ is an $\left( \frac{n}{4},2 \right)$-expander.
\end{enumerate}
\end{lemma}

\noindent {\bf Proof.} We will show that the subgraph $\Gamma _0$ constructed in Lemma \ref{sparse-P6prime} is an $\left( \frac{n}{4},2 \right)$-expander. Since is has at most $d_0n$ edges, this will finish the proof.\\
Let $U\subseteq V(G),\ |U|=k\leq \frac{n}{4}$. We will show that indeed $|N_{\Gamma _0}(U)| \geq 2k$.\\
Denote $U_1 = U\cap SMALL(G)$, $U_2 = U \setminus U_1$ and let $k_1,k_2$ be the sizes of $U_1,U_2$, respectively. Consider two cases:
\begin{enumerate}
\item $\frac{n}{\sqrt{\log n}}\leq k_2 \leq \frac{n}{4}$. We know that $|N_{\Gamma _0}(U_2)|\geq n-\frac{n}{\sqrt{\log n}}$, since otherwise the sets $U_2,\ V(G)\setminus N(U_2)$ contradict the conclusion of Lemma \ref{sparse-P6prime}. By  {\bf(P2)} we know that $k_1 \leq |SMALL(G)| \leq n^{0.3}$, so overall:
$$
|N_{\Gamma _0}(U)| \geq |N_{\Gamma _0}(U_2)\setminus U_1| \geq n-\frac{n}{\sqrt{\log n}}-n^{0.3} \geq \frac{n}{2} \geq 2k.
$$

\item $k_2 \leq \frac{n}{\sqrt{\log n}}$. By Lemma \ref{sparse-P6prime} we know that $\delta (\Gamma_0)\geq 2$, and by {\bf(P3)} we know that no two vertices in $U_1$ have distance less than 5. More specifically, every vertex of $U_1$ has at least two neighbours, and no two vertices share any common neighbour, so $|N_{\Gamma _0}(U_1)| \geq 2k_1$.\\
Next, we know that all vertices in $U_2$ have degree at least $d_0$ in $\Gamma _0$, and by {\bf(P4)} we have $e_{\Gamma _0}(U_2)\leq e_G(U_2) \leq k_2\cdot \log ^{\frac{3}{4}} n$, so:
$$
e(U_2,N_{\Gamma _0}(U_2)) \geq k_2\cdot d_0 - k_2\cdot \log ^{\frac{3}{4}} n \geq \frac{d_0k_2}{2}.
$$
This means that $|N_{\Gamma _0}(U_2)| \geq k_2 \cdot \log ^{\frac{1}{4}}n$, since otherwise $U_2,N_{\Gamma _0}(U_2)$ contradict {\bf(P5)}.\\
Finally, observe that, since the vertices of $U_1$ are at distance at least 5 from each other, we have:
$$
\forall u \in U_2:\ |(U_1\cup N_{\Gamma _0}(U_1))\cap (u \cup N_{\Gamma _0}(u))| \leq 1,
$$
which means $|(U_1\cup N_{\Gamma _0}(U_1))\cap (U_2 \cup N_{\Gamma _0}(U_2))| \leq k_2$. So overall:
\begin{eqnarray*}
|N_{\Gamma _0}(U)| & = & |N_{\Gamma _0}(U_2) \setminus U_1| + |N_{\Gamma _0}(U_1) \setminus (U_2 \cup N_{\Gamma _0}(U_2)|\\
& \geq & k_2\cdot \log ^{\frac{3}{4}}n - k_2 + 2k_1 - k_2 \geq 2k_1 + 2k_2 = 2k.
\end{eqnarray*}
\end{enumerate}
\hfill$\Box$

\begin{corollary}\label{sparse-corol}
Let $G$ be a graph such that {\bf(P0)}--{\bf(P7)} hold.  Then $G$ is Hamiltonian.
\end{corollary}

\noindent {\bf Proof.} By Lemma \ref{sparse-expander}, $G$ contains a subgraph $\Gamma _0$ which is an $\left( \frac{n}{4},2 \right)$-expander with at most $d_0n$ edges. By {\bf(P7)}, if $\Gamma _0$ is not Hamiltonian then it has a booster in $E(G) \setminus E(\Gamma _0)$. Add one such booster to $\Gamma _0$ to obtain a new subgraph of $G$, $\Gamma _1$. This new subgraph is also an  $\left( \frac{n}{4},2 \right)$-expander, this time with at most $d_0n +1$ edges, and so is either Hamiltonian or has a booster in  $E(G) \setminus E(\Gamma _1)$. This process can be repeated until we find a subgraph $\Gamma _i$ which is Hamiltonian. Since in each step the length of a longest path grows by at least $1$, at most $n$ steps are required, so $e(\Gamma _i)= e(\Gamma _0) + i \leq d_0n+n$, which means the process can always continue until Hamiltonicity is achieved.
\hfill$\Box$

\bigskip

Theorem \ref{main-thm} for $p \leq \frac{100\log n}{n}$ is obtained directly from Lemma \ref{sparse-properties} and Corollary \ref{sparse-corol}.

\subsection{The dense case}\label{subsec-dense}

Recall that we assume here that $\frac{100\log n}{n} \leq p \leq 0.01$. Define
$$
t_0 = 0.002\cdot np,
$$
and for a graph $G=(V,E)$, denote
$$
SMALL(G) = \lbrace v\in V(G):\, d(v) \leq t_0 \rbrace.
$$

Now, define the following graph properties:

\begin{itemize}
\item[{\bf (Q0)}] $\delta (G) \geq 2$;
\item[{\bf (Q1)}] $\Delta (G)\leq 7np$;
\item[{\bf (Q2)}] $|SMALL(G)|\leq 1$;
\item[{\bf (Q3)}] $\forall U,W\subset V(G)\ disjoint\ s.t.\ |U|=|W|= 10^{-13}n : e(U,W)\geq 10^{-27}{n^2p}$;
\item[{\bf (Q4)}] $\forall U,W\subset V(G)\ disjoint\ s.t.\ \frac{t_0}{3}-1 \leq |U| \leq 10^{-13}n,\ |W|=2|U|+3 : e(U,W) + e(U) < t_0|U|$;
\item[{\bf (Q5)}] Every $\left( \frac{n}{4},2 \right) $-expanding subgraph $\Gamma$ of $G$ with at most $(t_0 + 1)n$ edges is either Hamiltonian, or has a booster in $E(G) \setminus E(\Gamma )$.
\end{itemize}

\begin{lemma}\label{dense-properties}
Let $p=p(n)$ be such that $\frac{100\log n}{n} \leq p \leq 0.01$, and Let $G \sim G(n,p)$ be a random graph. Then the probability that all properties {\bf (Q0)}--{\bf (Q5)} hold is $1-(1+o(1))Pr(G \notin {\bf{(Q0)}})$.
\end{lemma}

\noindent {\bf Proof.} First, bound $Pr(G \notin {\bf{(Q0)}})$ from bellow. Recall from the proof of Lemma \ref{sparse-properties} that
$$
Pr(G\notin {\bf (Q0)}) \geq \sum_{v\in V(G)}Pr(A_v)\ -\sum_{u,v\in V(G)} Pr(A_u \cap A_v).
$$
Let $u,v \in V(G)$, so:
\begin{eqnarray*}
Pr(A_v) & = & Pr(d(v)=1)\\
& \geq & e^{-\frac{p}{1-p}(n-2)}np\\
& \geq & e^{-(1.1+o(1))np}\,;
\end{eqnarray*}
\begin{eqnarray*}
Pr(A_u \cap A_v) & \leq & Pr(e_G(\lbrace u,v \rbrace , V\setminus \lbrace u,v \rbrace ) \leq 2)\\
& \leq & e^{-(2-o(1))np}.
\end{eqnarray*}
And overall we get:
$$
Pr(G \notin {\bf{(Q0)}}) \geq n\cdot e^{-(1.1+o(1))np}-\binom{n}{2}e^{-(2-o(1))np} \geq e^{-(1.1+o(1))np}-e^{-1.5np} \geq e^{-1.2np}.
$$
Next, bound from above the probabilities of {\bf (Q1)}--{\bf (Q5)}:
\begin{itemize}
\item[{\bf (Q1).}] By Chernoff's inequality (Lemma \ref{chernoff}) with $\delta = 6$ and the union bound we get
$$
Pr(G \notin {\bf (Q1)}) \leq n\cdot Pr(Bin(n-1,p) \geq 7np) \leq ne^{-2np} = o(Pr(G \notin {\bf{(Q0)}}).
$$

\item[{\bf (Q2).}] The probability that $SMALL(G)$ contains more than one vertex is at most the probability that exist $u,v \in V(G)$ such that $e(\lbrace u,v \rbrace ,V(G)\setminus \lbrace u,v \rbrace ) \leq 2t_0$. So by the union bound:
\begin{eqnarray*}
Pr(G \notin {\bf (Q2)}) & \leq & \binom{n}{2} \left( 1 + \sum_{k=1}^{2t_0}\left( \frac{2enp}{k(1-p)} \right)^k \right) e^{-2(n-2)p}\\
& \leq & n^2t_0\cdot \left( \frac{enp}{0.99t_0} \right)^{2t_0} e^{-2(n-2)p}\\
& \leq & \exp (3\log n + 15t_0-2np+2p) \leq e^{-1.9np},
\end{eqnarray*}
which is sufficiently small.

\item[{\bf (Q3).}] By Chernoff's inequality and the union bound we get that  $Pr(G \notin {\bf (Q3)})$ is at most
$$
\binom{n}{10^{-13}n}^2\cdot Pr\left( Bin\left( 10^{-26}n^2,p \right) \leq 10^{-27}n^2p\right) \leq 4^ne^{-\frac{1}{4}n^2p} = e^{-\Omega (n^2p)},
$$
and thus is small enough.

\item[{\bf (Q4).}] By the union bound:
\begin{eqnarray*}
Pr(G \notin {\bf (Q4)}) & \leq & \sum_{k=\frac{t_0}{3}-1}^{10^{-13}n} \binom{n}{k}\binom{n}{2k+3} \left( \frac{ep\cdot \left( k(2k+3) + \binom{k}{2} \right) }{t_0k} \right)^{t_0k}\\
& \leq & \sum_{k=\frac{t_0}{3}-1}^{10^{-13}n} \left( n^4 \cdot e^{-21t_0} \right) ^k  \leq  \sum_{k=\frac{t_0}{3}-1}^{10^{-13}n} e^{-t_0k}\\
& \leq & (1+o(1))e^{-\Omega ({t_0}^2)} = o(Pr(G \notin {\bf{(Q0)}}).
\end{eqnarray*}

\item[{\bf (Q5).}] Similarly to {\bf (P7)} in the proof of Lemma \ref{sparse-properties}, the probability of the existence of such a subgraph $\Gamma\subseteq G$, with none of its boosters being an edge of $G$, is at most
\begin{eqnarray*}
Pr(G \notin {\bf (P7)}) & \leq & \sum_{k=1}^{(t_0+1)n} \binom{\binom{n}{2}}{k}p^k (1-p)^{\frac{n^2}{32}}\\
& \leq & e^{o(n)}\cdot \left( \frac{enp}{2t_0} \right)^{(t_0+1)n} e^{-\frac{n^2p}{32}} \\
& \leq & e^{o(n)}\cdot \left( \frac{e}{0.004} \right)^{0.003n^2p} e^{-\frac{n^2p}{32}}\\
& \leq & \exp ((0.025 - 0.03)n^2p) = o(Pr(G\notin {\bf (P0)})).
\end{eqnarray*}
\end{itemize}
\hfill$\Box$

\begin{lemma}\label{dense-P3prime}
Let $G$ be a graph such that {\bf(Q0)}, {\bf(Q1)} and {\bf(Q3)} hold.  Then $G$ contains a subgraph $\Gamma_0$ such that:
\begin{enumerate}
\item $\delta (\Gamma_0) \geq 2$,
\item $\Gamma_0$ has at most $t_0n$ edges,
\item $\forall v\notin SMALL(G) \,:\, d_{\Gamma _0}(v) \geq t_0$,
\item $\forall U,W\subset V(G)\ disjoint\ vertex\ sets\ s.t.\ |U|=|W|= 10^{-13}n : e_{\Gamma _0}(U,W)\geq 1$.
\end{enumerate}
\end{lemma}

\noindent {\bf Proof.} Recall the random subgraph construction described in the proof of Lemma \ref{sparse-P6prime}. Consider the same construction, with $t_0$ replacing $d_0$. The probability of the existence of $U,W\subset V(G)$ with $|U|=|W|= 10^{-13}n $ such that no edges cross between them is at most
$$
\exp \left( -\frac{t_0}{\max_{u\in U}d_G(u)}\cdot e_G(U,W)+o(n) \right),
$$
and by assuming $G \in {\bf(Q1)},{\bf(Q3)}$ we finally get
$$
\leq \exp \left( -\frac{0.002}{7}10^{-13}n^2p + o(n) \right) = o(1).
$$
\hfill$\Box$

\begin{lemma}\label{dense-expander}
Let $G$ be a graph such that {\bf(Q0)}--{\bf(Q4)} hold.  Then $G$ contains a subgraph $\Gamma_0$ such that:
\begin{enumerate}
\item $\Gamma_0$ has at most $t_0n$ edges,
\item $\Gamma_0$ is an $\left( \frac{n}{4},2 \right)$-expander.
\end{enumerate}
\end{lemma}

\noindent {\bf Proof.} We will show that the subgraph $\Gamma_0$ constructed in Lemma \ref{dense-P3prime} is an $\left( \frac{n}{4},2 \right)$-expander.\\
Let $U \subset V(G)$ be a vertex subset of size $k \leq \frac{n}{4}$. We will show that $|N_{\Gamma_0}(U)| \geq 2k$. Consider the following cases:

\begin{enumerate}
\item $k=1$. Since $\delta (\Gamma_0) \geq 2$ we get $|N(U)| \geq 2 = 2k$.
\item $2 \leq k \leq \frac{t_0}{3}$. Recall that by the construction of $\Gamma_0$, the degree in $\Gamma_0$ of any vertex $v\in V(G)\setminus SMALL(G)$ is at least $t_0$. By {\bf(Q2)}, $|SMALL(G)| \leq 1$, so $U\setminus SMALL(G)$ is non-empty. Let $v$ be a member of $U\setminus SMALL(G)$. It holds that
$$
|N(U)| \geq |N(\lbrace v \rbrace ) \setminus e_{\Gamma_0}(v,U)| \geq t_0 - k \geq 2k.
$$
\item $\frac{t_0}{3} \leq k \leq 10^{-13}n$. Let $U^{\prime} := U \setminus SMALL(G)$ and $W:=N(U) \cup (U \cap SMALL(G))$. From {\bf(Q2)} we have $\frac{t_0}{3}-1 \leq |U^{\prime}| \leq 10^{-13}n$. Since all vertices of $U^{\prime}$ have degree at least $t_0$, and all edges touching $U^{\prime}$ have their other end in either $U^{\prime}$ or $W$, we get that $e(U^{\prime}) + e(U^{\prime},W) \geq t_0|U^{\prime}|$. By {\bf(Q4)} this means
$$
|N(U)| \geq |W|-1 \geq 2|U^{\prime}|+3-1 \geq 2k.
$$
\item $10^{-13}n \leq k \leq \frac{n}{4}$. Observe that by Lemma \ref{dense-P3prime}: $|V(G) \setminus (U \cup N(U))| \leq 10^{-13}n$. So
$$
|N(U)| \geq n-k-10^{-13} \geq \frac{n}{2} \geq 2k.
$$
\end{enumerate}
\hfill$\Box$

\begin{corollary}\label{dense-corol}
Let $G$ be a graph such that {\bf(Q0)}--{\bf(Q5)} hold.  Then $G$ is Hamiltonian.
\end{corollary}

\noindent {\bf Proof.} The proof is essentially identical to the proof of Corollary \ref{sparse-corol}.\hfill$\Box$

\bigskip

Theorem \ref{main-thm} for $\frac{100\log n}{n}\leq p \leq 0.01$ is obtained from Lemma \ref{dense-properties} and Corollary \ref{dense-corol}.

\subsection{The very dense case}\label{subsec-vdense}

Recall that we assume here that $p \geq 0.01$. Define the following graph properties:

\begin{itemize}
\item[{\bf (R0)}] $\delta (G) \geq 2$;
\item[{\bf (R1)}] There is at most one vertex $v \in G$ s.t. $d(v) \leq \frac{np}{10}$;
\item[{\bf (R2)}] $\forall U,W\subset V(G)\ disjoint\ s.t.\ |U|=|W|= \frac{np}{30} : e(U,W)\geq 1$;
\item[{\bf (R3)}] $\alpha (G) \leq \frac{np}{40}$.
\end{itemize}

\begin{lemma}\label{vdense-properties}
Let $p=p(n)$ be such that $p \geq 0.01$, and Let $G \sim G(n,p)$ be a random graph. Then the probability that all properties {\bf (R0)}--{\bf (R3)} hold is $1-(1+o(1))Pr(G \notin {\bf{(R0)}})$.
\end{lemma}

\noindent {\bf Proof.} First, bound $Pr(G \notin {\bf{(R0)}})$ from below by Bonferroni's inequality.\\
Let $u,v\in V(G)$, so:
\begin{eqnarray*}
Pr(A_v) & = & Pr(d(v)=0)+Pr(d(v)=1) \\
& = & (1-p)^{n-1} + (n-1)p(1-p)^{n-2} \\
& = & (1-p)^{(1+o(1))n}\,;
\end{eqnarray*}
\begin{eqnarray*}
Pr(A_u \cap A_v) & \leq & Pr(e_G(\lbrace u,v \rbrace , V\setminus \lbrace u,v \rbrace ) \leq 2)\\
& = & (1-p)^{(2-o(1))n}.
\end{eqnarray*}
And overall we get:
$$
Pr(G \notin {\bf{(R0)}}) \geq n\cdot (1-p)^{(1+o(1))n}-\binom{n}{2}(1-p)^{(2-o(1))n} = (1-p)^{(1+o(1))n}.
$$
Next, bound from above the probabilities of {\bf (R1)}--{\bf (R3)}:
\begin{itemize}
\item[{\bf (R1).}] The probability that $G$ contains more than one vertex of degree at most $\frac{np}{10}$ is at most the probability of the existence of $u,v \in V(G)$ such that $e(\lbrace u,v \rbrace ,V(G)\setminus \lbrace u,v \rbrace ) \leq \frac{np}{5}$. So by the union bound:
\begin{eqnarray*}
Pr(G \notin {\bf (R1)}) & \leq & \binom{n}{2} \left( 1 + \sum_{k=1}^{0.2np}\left( \frac{2enp}{k(1-p)} \right)^k \right) (1-p)^{2(n-2)}\\
& \leq & ( 10e )^{0.2np} (1-p)^{(1.8-o(1))n}\\
& \leq & (1-p)^{(1.1-o(1))n}\cdot \exp \left( (0.2\log (10e) - 0.7)np \right)\\
& \leq & (1-p)^{(1.1-o(1))n} = o(Pr(G \notin {\bf{(R0)}})).
\end{eqnarray*}
\item[{\bf (R2).}] By the union bound:
$$
Pr(G \notin {\bf (R2)}) \leq \binom{n}{\frac{np}{30}}^2 (1-p)^{ \left( \frac{np}{30} \right) ^2} = (1-p)^{\Omega \left( n^2 \right)} = o(Pr(G \notin {\bf{(R0)}})).
$$
\item[{\bf (R3).}] Similarly, by the union bound:
$$
Pr(G \notin {\bf (R3)}) \leq \binom{n}{ \frac{np}{40} } (1-p)^{\binom{\frac{np}{40}}{2}} = (1-p)^{\Omega \left( n^2 \right)} = o(Pr(G \notin {\bf{(R0)}})).
$$
\hfill$\Box$

\end{itemize}

\begin{lemma}\label{vdense-ham}
Let $G$ be a graph such that {\bf(R0)}--{\bf(R3)} hold.  Then $G$ is Hamiltonian.
\end{lemma}

\noindent {\bf Proof.} Denote by $v$ a vertex of $G$ such that $d(v) = \delta (G)$, and by $G^{\prime}$ the graph obtained from $G$ by removing $v$ and all its edges.\\
We observe that the properties {\bf(R2)}, {\bf(R3)} hold for $G^{\prime}$, since removing vertices does not affect these properties. Furthermore, by property {\bf(R1)} we have $\delta (G^{\prime}) \geq \frac{np}{10}-1$. We claim that these are sufficient for showing that $\kappa (G^{\prime}) \geq \frac{np}{30}$.\\
Indeed, suppose towards contradiction that there exists a set $U\subseteq V(G^{\prime})$ of size $\frac{np}{30}$ such that removing the vertices of $U$ disconnects $G^{\prime}$, and denote by $W_1, W_2$ two components in the resulting graph (WLOG we assume $|W_1|\leq |W_2|$). Consider the following cases:
\begin{enumerate}
\item $|W_1| \leq \frac{np}{30}$. Let $w\in W_1$. Since $W_1$ is a component of  $G^{\prime} - U$, we know that $N_{G^{\prime}}(w) \subseteq W_1\cup U$. But
$$
|N_{G^{\prime}}(w)| = d_{G^{\prime}}(w) \geq \frac{np}{10} - 1 > |W_1|+|U|\,,
$$
a contradiction.
\item $|W_1| > \frac{np}{30}$. Since $|W_2|\geq |W_1| > \frac{np}{30}$ and $e_{G^{\prime }}(W_1,W_2)=0$, this is a contradiction to {\bf(R2)}.
\end{enumerate}
Now, from {\bf(R3)} we get that $\alpha (G^{\prime}) \leq \frac{np}{40} < \kappa (G^{\prime})$. By Theorem \ref{chvatal-erdos} this means that $G^{\prime}$ is Hamilton--connected.\\
We now return to $G$. By {\bf(R0)}, $v$ has at least two neighbours, say $u_1,u_2$. Since $G^{\prime}$ is Hamilton--connected, it contains a Hamilton path $\mathcal{P}$ with $u_1,u_2$ being its two endpoints. Now $(u_1,v),(v,u_2),\mathcal{P}$ is a Hamilton cycle in $G$.\hfill$\Box$

\section{Perfect matching}\label{sec-PM}

We provide a sketch of a proof for Theorem \ref{matching-thm}.\\
Observe that similarly to the case of Hamilton cycles, for the very sparse case, here defined as $p(n)$ such that $np-\log n$ does not tend to infinity, the result is already known. This is due to the classical result by Erd\H{o}s and R{\'e}nyi, stating:

\[
\lim _{n \rightarrow \infty }\Pr(G\ \mbox{contains a perfect matching})=
\begin{cases}
    0 & np - \log n \rightarrow -\infty,\\
    e^{-e^{-c}} & np - \log n \rightarrow c,\\
    1 & np - \log n \rightarrow \infty.
\end{cases}
\]
Which means that $\lim _{n \rightarrow \infty }\Pr(G\ \mbox{contains no perfect matching}) = \lim _{n \rightarrow \infty }\Pr(\delta (G) =0) = C > 0$, which suffices.\\
We also observe that the very dense case of $p \geq 0.01$ was proven in Section \ref{subsec-vdense} up to a small adjustment. Recall that in the proof of Lemma \ref{vdense-ham} we showed that in a graph with properties {\bf(R1)}--{\bf(R3)}, removing a minimum degree vertex yields a Hamilton--connected graph. Replace the property {\bf(R0)} with a new property {\bf (R0')}: $\delta (G) > 0$. By taking the minimum degree vertex and a Hamilton path starting at one of it's neighbours, we get a Hamilton path in $G$, which becomes a perfect matching by taking every other edge. Since like {\bf(R0)}, $Pr(G \notin {\bf (R0\mbox{\bf{'}})}) = (1-p)^{(1+o(1))n}$, we obtain the desired result.\\
Towards the goal of providing a proof idea for other ranges of $p(n)$, we introduce the notion of staples, which will be used here similarly to the way boosters were used in the proof of Theorem \ref{main-thm}:

\begin{defin}\label{def-staple}
Let $G=(V,E)$ be a graph. A non-edge $(u,v)\in E(G)$ is called a {\em staple} if the graph $G^{\prime}$ with edge set $E(G^{\prime})=E(G)\cup \lbrace (u,v) \rbrace$ has a perfect matching or a matching larger than a maximum matching of $G$.
\end{defin}

\noindent Similarly to Lemma \ref{lemma-booster}, in \cite{FK} the following is proved (see Lemma 6.3):

\begin{lemma}\label{lemma-staple}
Let $G=(V,E)$ be a graph with no perfect matching, and assume that $G$ is a $(k,1)$-expander. Then $G$ has at least $\binom{k+1}{2}$ staples.
\end{lemma}

\noindent This enables us to take an approach similar to the proof in Section \ref{sec-proof}, following similar steps:

\begin{itemize}
\item[(i)] Show that for $G \sim G(n,p)$, the probability that $G$ does not contain a relatively sparse $(k,1)$--expanding subgraph, for some $k$ linear in $n$, is $(1+o(1))Pr(\delta (G) = 0)$;
\item[(ii)] Show that the probability that $G$ does not contain a staple for each of its sparse $(k,1)$--expanding subgraphs is $o(Pr(\delta (G) = 0))$.
\end{itemize}
Recalling Corollary \ref{sparse-corol}, this suffices for proving Theorem \ref{matching-thm}.


\begin{thebibliography}{99}

\bibitem{AKS85} M. Ajtai, J. Koml{\'o}s and E. Szemer{\'e}di, \textit{First occurrence of Hamilton cycles in random graphs}, Cycles in graphs '82, North Holland Mathematical Studies 115, North Holland, Amsterdam (1985), 173--178.

\bibitem{B84} B. Bollob{\'a}s, \textit{The evolution of sparse graphs}, Graph Theory and Combinatorics, Academic Press, London (1984), 35--57.

\bibitem{BT} B. Bollob{\'a}s and A. Thomason, \textit{Random graphs of small order}, Random graphs '83 (Pozna\'n, 1983), North Holland Mathematical Studies 118, North Holland, Amsterdam (1985), 47--97.

\bibitem{CE}  V. Chv{\'a}tal and P. Erd\H{o}s, \textit{A note on Hamiltonian circuits}, Discrete Mathematics 2 (1972), 111--113.

\bibitem{ER} P. Erd\H{o}s and A. R{\'e}nyi, \textit{On the existence of a factor of degree one of a connected random graph}, Acta Mathematica Academiae Scientiarum Hungaricae 17 (1966), 359--368.

\bibitem{FK} A. Frieze and M. Karo{\'n}ski, {\bf Introduction to random graphs}.  Cambridge University Press, 2015.

\bibitem{HKS09} D. Hefetz, M. Krivelevich and T. Szab{\'o}. \textit{Hamilton cycles in highly connected and expanding graphs}, Combinatorica 29 (2009), 547--568.

\bibitem{CHER} W. Hoeffding, \textit{Probability inequalities for sums of bounded random variables}, Journal of the American Statistical Association 58 (1963), 13--30.

\bibitem{KS83} J. Koml{\'o}s and E. Szemer{\'e}di, \textit{Limit distributions for the existence of Hamilton circuits in a random graph}, Discrete Mathematics 43 (1983), 55--63.

\bibitem{K} M. Krivelevich, \textit{Long paths and Hamiltonicity in random graphs}, Random Graphs, Geometry and Asymptotic Structure, London Mathematical Society Student Texts 84, Cambridge University Press (2016), 4--27.

\bibitem{MY16} C. McDiarmid and N. Yolov, \textit{Hamilton cycles, minimum degree and bipartite holes}, Journal of Graph Theory 86 (2017), 277--285.

\bibitem{POS} L. P{\'o}sa, \textit{Hamiltonian circuits in random graphs}, Discrete Mathematics 14 (1976), 359--364.

\end{thebibliography}
\end{document}